\documentclass
{notices}
\usepackage[utf8]{inputenc}
\usepackage[T1]{fontenc}
\usepackage{amssymb,amscd, amsthm, amsmath}
\usepackage[all,line,arc,curve,color,frame,pdf]{xy}
\usepackage{tikz}
\usepackage{tikz-cd}
\usetikzlibrary{positioning, trees, snakes}
\usepackage[textsize=tiny]{todonotes}
\usepackage{hyperref}
\usepackage[shortlabels]{enumitem}
\usepackage{float} 
\usepackage[paper=a4paper, margin=3.5cm]{geometry}
\usepackage{cleveref}   
\usepackage{comment}
\usepackage{scalerel}
\usepackage[normalem]{ulem}
\usepackage{listings}
\usepackage{verbatim}
\usepackage{marvosym} 

\DeclareMathOperator{\Pic}{Pic}
\DeclareMathOperator{\CH}{CH}

\def\ZZ{\mathbb{Z}}
\def\CC{\mathbb{C}}
\def\PP{\mathbb{P}}
\def\QQ{\mathbb{Q}}
\def\RR{\mathbb{R}}
\theoremstyle{theorem}
\newtheorem{theorem}{Theorem}[section]
\theoremstyle{definition}
\newtheorem{example}{Example}[section]
\newtheorem{definition}{Definition}[section]
\newtheorem{conjecture}{Conjecture}[section]

\begin{document}
\title{Enumerative geometry meets statistics, combinatorics and topology}

\author{Mateusz Micha{\l}ek}
  \maketitle
\affil{
	University of Konstanz, Germany, Fachbereich Mathematik und Statistik, Fach D 197
	D-78457 Konstanz, Germany
}
\thanks{The author is funded by the Deutsche Forschungsgemeinschaft –- Projektnummer 467575307.}



\begin{abstract} 
We explain connections among several, a priori unrelated, areas of mathematics: combinatorics, algebraic statistics, topology, and enumerative algebraic geometry. Our focus is on discrete invariants, strongly related to the theory of Lorentzian polynomials.
The main concept joining the mentioned fields is a linear space of matrices.    
\end{abstract}


\section{Introduction}
The following questions arise prominently in different branches of mathematics.
\begin{enumerate}
\item Given a graph $G$, how many proper vertex colorings with $k$ colors exist?
\item What is the degree of a general linear concentration model? What is its maximum likelihood degree?
\item What is the Euler characteristic of a hypersurface defined by the determinant of a matrix with linear entries?
\item How many degree two hypersurfaces pass through $a$ general points and are tangent to $b$ general hyperplanes?
\end{enumerate}
We will start by explaining the meaning of the above questions. 
Our main aim is to show that the central objects we encounter are in fact shadows of one construction and all of the above questions are in fact one (or more precisely two related) question(s). It turns out that the unifying setting is surprisingly simple: we will always start from a linear space $L$ of square matrices. To such a space we canonically associate a polynomial with integral coefficients. This is a special instance of the so-called volume polynomial and an example of a \emph{Lorentzian polynomial}. Roughly speaking, the coefficients of this polynomial form the multidegree of the graph of the gradient of the determinant restricted to $L$. We also present powerful, modern geometric tools to study these basic invariants from a new perspective. This will be achieved by performing intersection theory on smooth, projective, and beautiful varieties. 
We will work over the field $\CC$.
\section*{Acknowledgements}
We thank Rodica Dinu, Christopher Eur, Maciej Ga{\l}{\k a}zka, Lukas Gustafsson, Kangjin Han, Joachim Jelisiejew, Anna Seigal, Tim Seynnaeve, Evgeny Shinder and Emanuele Ventura for important comments on the first version of this article. 

\section{Main players}
\subsection*{Multidegree}
For our purposes a variety $V(f_1,\dots,f_n)$ is the set of zeros of a polynomial system $f_i=0$. Given homogeneous polynomials on a vector space $V$, it is natural to consider their zero set as a subset of the projective space $\PP(V)$. The corresponding variety is called \emph{a projective variety} \cite[Chapters 1-3]{michalek2021invitation}.

\begin{example}\label{ex:quadric}
Let $f=xy-zt$. We may regard $V(f)$ as a $2$-dimensional quadratic surface in $\PP^3$. 
\end{example}

Citing Bernd Sturmfels ``\emph{The two most important invariants of a variety in a projective space are its dimension and degree}''. One of a few equivalent ways to define the degree of a variety $X\subset \PP^n$ is as the number of points one obtains after intersecting it with $\dim X$-many general hyperplanes. A variety $X$ is \emph{irreducible} if, whenever it is a union of two varieties $X=X_1\cup X_2$ then $X=X_1$ or $X=X_2$. To the projective space $\PP^n$ one associates a ring $H^*(\PP^n)$, which in this case coincides with the Chow ring and the cohomology ring. Elements of this ring are formal linear combinations of classes of irreducible subvarieties of $\PP^n$, called cohomology classes. A variety is equivalent to a formal sum of its irreducible components. Further, two, possibly reducible, varieties are equivalent, if they have the same degree and all components are of the same dimension. The multiplication in $H^*(\PP^n)$ corresponds to intersection of varieties\footnote{Assuming we choose representatives of classes that intersect in a nice way --- formally we have to assume they intersect transversally.}.
Note that all hyperplanes define the same cohomology class. In fact, the cohomology ring with rational coefficients of $\PP^n$ is $\QQ[H]/(H^{n+1})$, where $H$ is the class of the hyperplane. Then the class of an irreducible variety $X$ equals $(\deg X)H^{n-\dim X}$. 

\begin{example}
Continuing the example of $f=xy-zt$, one can check that $V(f)$ has degree two. In general, the notion of the degree of a variety generalizes the notion of the degree of a polynomial.
\end{example}

When $X$ is a subvariety of a product $\PP^n\times \PP^m$ of projective spaces, the analogue of the degree is the \emph{multidegree}. Indeed, we have two different (families of) ``hyperplanes'' in $\PP^n\times \PP^m$. Namely, the product $H_1$ of a hyperplane in $\PP^n$ with $\PP^m$, and  the product $H_2$ of $\PP^n$ with a hyperplane in $\PP^m$. Thus, instead of one number, we obtain a sequence of $(\dim X+1)$-many numbers, by intersecting $X$ with $a$ general hyperplanes of type $H_1$ and $b$ general hyperplanes of type $H_2$, where $a+b=\dim X$. Analogously to the previous case, the cohomology ring of $\PP^n\times \PP^m$ is $\QQ[H_1,H_2]/(H_1^{n+1},H_2^{m+1})$. The multidegree of $X$ tells us the cohomology class of $X$. For more information about the multidegree we refer to the book \cite[Chapter 8]{miller2005combinatorial}.

\subsection*{Graphs} Let $G=(V,E)$ be a loopless graph. A \emph{proper vertex coloring} using $k\in\ZZ_{\geq 0}$ colors is a function $f:V\rightarrow \{1.\dots,k\}$, such that whenever two vertices $v_1,v_2$ are connected by an edge, we have $f(v_1)\neq f(v_2)$. The function $\chi_G:\ZZ_{\geq 0}\rightarrow\ZZ_{\geq 0}$ that to $k$ assigns the number of proper vertex colorings using $k$ colors is known as the \emph{chromatic polynomial} of $G$. It is indeed a polynomial, as one may prove by induction on the number of edges, by contracting and deleting a given edge. The same proof shows that $\chi_G$ is a polynomial of degree $|V|$ with integral coefficients, whose signs are alternating. If $G$ has at least one edge then $\chi_G$ has a root at one. We thus define the \emph{reduced chromatic polynomial} $\overline{\chi}_G(k):=\chi_G(k)/(k-1)$. From now on, for simplicity, we will assume that $G$ is connected.
\begin{example}\label{exm:chiC3}
Let $P_n$ be the path with $n$ vertices. Given $k$ colors we may assign to the first vertex any of the colors. Then to each consecutive vertex we may assign $(k-1)$ colors. We obtain:
\begin{align*}
\chi_{P_n}(k)&=k(k-1)^{n-1}\\
\overline{\chi}_{P_n}(k)&=k(k-1)^{n-2}.
\end{align*}

Let $C_3$ be the $3$-cycle, i.e.~a triangle. If we remove one edge, we obtain a path $P_3$ with $3$ vertices. 
A proper coloring of $P_3$ is not a proper coloring of $C_3$ if and only if the two end vertices have the same color. Thus, there are exactly $\chi_{P_2}(k)=k(k-1)$ such non-proper colorings. We obtain:
\begin{align*}
\chi_{C_3}(k)&=\chi_{P_3}(k)-\chi_{P_2}(k)=k(k-1)(k-2)\\
\overline{\chi}_{C_3}(k)&=k(k-2).
\end{align*}
\end{example}
\begin{example}\label{exm:cycle}
Let $C_4$ be the $4$-cycle and let $P_4$ be the path obtained by removing one edge from $C_4$. A proper coloring of $P_4$ is not a proper coloring of $C_4$ if and only if the two end vertices have the same color. Thus, there are exactly $\chi_{C_3}(k)=k(k-1)(k-2)$ such colorings. We obtain:
\begin{align*}
\chi_{C_4}(k)&=k(k-1)\left( (k-1)^2-(k-2)\right)=\\
&k(k-1)(k^2-3k+3)\\
\overline{\chi}_{C_4}(k)&=k^3-3k^2+3k.
\end{align*}
\end{example}
Let $n:=|E|$. Below, we describe a classical construction of associating to $G$ a subspace of $\CC^n$, which is a special case of a representation of a matroid. Let $\CC^n$ be the vector space with basis $f_e$ for $e\in E$. We orient the edges of $G$ in an arbitrary way. Let $L_G$ be the subspace of $\CC^n$ spanned by all vectors $w_v$ indexed by vertices $v\in V$, where the coordinate of $w_v$ corresponding to the edge $e$ is given by: \[f_e^*(w_v)=\begin{cases}
1\text{ when }e=(v,v')\\
-1\text{ when }e=(v',v)\\
0\text{ otherwise.}
\end{cases}\]
\begin{example}
When $C_4$ is the oriented $4$-cycle we see that $L_G$ is a codimension one subspace of $\CC^4$ given by the linear equation: sum of coordinates equal to zero.
\end{example}
A magnificent connection of invariants of $L_G$ with the chromatic polynomial $\chi_G$ was discovered by Huh \cite{Huh1}. This was a cornerstone to solutions of several long-standing conjectures allowing applications of powerful theorems from algebraic geometry to combinatorics.  We describe this connection below using a language different from \cite{Huh1}, however better suited to exhibit connections with other topics. First, the ambient space $\CC^n$ will be identified with the space $M_{n}^D$ of $n\times n$ diagonal matrices. On the projective $(n-1)$-dimensional space $\PP(M_n^D)$ we have a rational\footnote{i.e.~a map that is not defined everywhere} map:
\[F:\PP(M_n^D)\dashrightarrow \PP^{n-1}\]
given by the (nonzero) $(n-1)\times(n-1)$ minors. At this point the reader should see that this is simply the classical Cremona transformation that inverts coordinates. One could also say this is the gradient of the determinant. We restrict $F$ to $\PP(L_G)$
and look at the graph:
\[\Gamma_G:=\overline{\{(x,y)\in \PP(L_G)\times\PP^{n-1}: y=F(x)\}}.\]
The dimension of $\Gamma_G$ is the dimension of $\PP(L_G)$. For connected graphs it simply equals the cardinality of the vertex set minus two\footnote{Experts may recognize that this is in fact the rank of the matroid minus one.}. The multidegree turns out to be utterly interesting!
\begin{theorem}[\cite{Huh1}]
The multidegree sequence of $\Gamma_G$ equals the sequence of absolute values of coefficients of the reduced chromatic polynomial $\overline{\chi}_G$.
\end{theorem}
\begin{example}\label{MD_C4}
Continuing the example when $C_4$ is the $4$-cycle we obtain the restriction of the map $\PP^3\dashrightarrow \PP^3$:
\[(x:y:z:t)\mapsto (yzt:xzt:xyt:xyz)\]
to $\PP^2\subset\PP^3$ given by $x+y+z+t=0$. The image of the restriction is the hypersurface defined by the degree three elementary symmetric polynomial. This degree three is also one entry of the multidegree; it is given by the product $[\Gamma_{C_4}] [H_2]^{2}=3[pt]$ in the cohomology ring $H^*(\PP^3\times\PP^3)$, where $[pt]$ is the class of a point. It also corresponds to the lowest degree term $3k$ in $\overline{\chi_{C_4}}$ in Example \ref{exm:cycle}. The computation of the multidegree can be achieved e.g.~through the \texttt{Cremona} package in Macaulay2 \cite{M2}.
\begin{verbatim}
R=QQ[a,b,c,d], S=QQ[x,y,z,t]
f=map(S/ideal(x+y+z+t),R,
{y*z*t,x*z*t,x*y*t,x*y*z})
loadPackage"Cremona"
projectiveDegrees f
\end{verbatim}
The output, consistent with Example \ref{exm:cycle}, is:
\begin{verbatim}
o4 = {1,3,3}
\end{verbatim}
\end{example}

\subsection*{Statistical models}
We start with the following thought experiment. We are given a coin with probability of heads $p$ and tails $1-p$. Say we know that $p$ was chosen from the interval $(\frac{1}{3},\frac{2}{3})$. Such an assumption corresponds to choosing a \emph{statistical model}. How to find the value of $p$? We could start throwing the coin many times. Say, after throwing $1000$ times we had $400$ heads. Intuitively we estimate $p=0.4$. A rigorous way is to compute the probability of our experiment --- i.e.~$400$ heads --- as a function of $p$. This is referred to as the likelihood function, which is $L(p):=p^{400}(1-p)^{600}$. The $p$ we want to find \emph{maximizes the likelihood function}. Indeed, it may be checked that $p=0.4$ is the correct one. This method is called \emph{maximum likelihood estimation}. In practice, one often maximizes the log-likelihood function $\log L(p)$. An efficient way is to check when the derivative is equal to zero. In the described case:
\[\frac{d(\log L)}{dp}=\frac{200(2-5p)}{p(1-p)}\]
equals zero precisely when $p=0.4$. For more complicated statistical models we may get several critical values, one of which is the desired maximum (assuming it is achieved). The number of complex critical points is known as the \emph{maximum likelihood degree} (ML-degree) and is one of the basic algebraic measures of the complexity of the model \cite{StUh}. 

Our main interest will be \emph{Gaussian models}. As we will soon see the degree and the ML-degree of linear multivariate Gaussian models will be governed by the multidegree of a graph of a map given by the gradient of the determinant. The classical Gaussian distribution on $\RR$ has density given by:
\[f(x)=\frac{1}{\sigma\sqrt{2\pi}} e^{-\frac{1}{2}\left(\frac{x - \mu}{\sigma}\right)^2}.\]
Here $\mu\in \RR$ is the mean and $\sigma^2\in\RR_+$ is the variance. A slight generalization is a multivariate Gaussian distribution, i.e.~distribution on $\RR^n$. In particular, the mean $\mu\in\RR^n$ is now a vector. How should we generalize the variance? We note that the exponent in $f(x)$ is actually an evaluation of a quadratic form on $x-\mu$. In the one-dimensional case we only had to choose one coefficient for the quadratic form. For $\RR^n$ the right generalization is a positive-definite symmetric $n\times n$ matrix $\Sigma$, known as the covariance matrix. We obtain:
\begin{align*}&f_{\Sigma,\mu}(x)=\\
&\det(2\pi\boldsymbol\Sigma)^{-\frac{1}{2}} \, \exp \left( -\frac{1}{2}(\mathbf{x} - \boldsymbol\mu)^{{{\!\mathsf{T}}}} \boldsymbol\Sigma^{-1}(\mathbf{x} - \boldsymbol\mu) \right).
\end{align*}
  As $\Sigma$ is positive-definite the function is integrable. The coefficient $\det(2\pi\boldsymbol\Sigma)^{-\frac{1}{2}}$ in front of the exponential is chosen so that the integral of $f_{\Sigma,\mu}(x)$ over $\RR^n$ equals one. Thus indeed we obtain a probability distribution.

In the examples above, we were not specifying one probability distribution, but a set of those. This motivates the following definition.
\begin{definition} A \emph{statistical model} is a family of probability distributions.
\end{definition}
Such a family often comes with an additional structure, e.g.~as a subset of $\RR^m$. For the general multivariate Gaussian model we would take $m=\binom{n+1}{2}+n$ and identify the probability distribution with a point $(\Sigma,\mu)\in \RR^m$. Even when dealing with Gaussian distributions we often make further assumptions on $\Sigma$ and $\mu$. First, we note that if we have data that we believe come from a multivariate Gaussian distribution, we may estimate the mean $\mu$ by taking the mean of the data. By shifting the data we will assume from now on that $\mu=0\in\RR^n$. Thus, the model is specified by determining a set $S\subset \RR^{\binom{n+1}{2}}$ to which $\Sigma$ may belong. Of particular interest for us will be the \emph{linear concentration model} described below. 
The matrix $K:=\Sigma^{-1}$ is called the \emph{concentration matrix}. Let us fix a linear space $L$ of symmetric $n\times n$ matrices that contains a positive-definite matrix. A linear concentration model is given by:
\[S_L:=\{\Sigma: K=\Sigma^{-1}\in L\}.\]
Formally we still require $\Sigma$, and hence also $K$, to be positive definite.

We have two important invariants of $S_L$. One is the \emph{degree of the model}, which is simply the degree of the variety $L^{-1}:=\overline{S_L}$ that is the Zariski closure of $S_L$ in the ambient space of $n\times n$ matrices. We note that $\overline{S_L}$ is also the Zariski or Euclidean closure of the locus of inverses of all invertible  matrices in $L$. The other invariant is the ML-degree that we define below in analogy to the case of a coin.   Suppose we are given a data vector $d_1\in \RR^n$. It does not make sense to ask what is the probability of observing $d_1$, as this is zero. Still, the value of the density function on $d_1$ is the correct measure of how likely it is to make such an observation.
Hence, for general data $d_1,\dots,d_k\in \RR^n$ with mean zero the \emph{ML-degree} is the number of (complex) critical points of the log-likelihood function:
\[\Sigma\mapsto \log\left(\prod_{i=1}^k f_{\Sigma,0}(d_i)\right).\]
As our probability distribution is no longer discrete, as it was in case of the coin, instead of maximizing the probability we maximize the product of values of density functions. If no condition on $\Sigma$ is required, i.e.~$L$ is the whole ambient space, one can check \cite[p.43 and 44]{drton2008lectures} that the optimal $\Sigma$ is given by:
\[\hat\Sigma:=\frac{1}{k} \sum_{i=1}^k d_i d_i^t\]
and is called \emph{the sample covariance matrix}. 

Which $\Sigma\in S_L$ maximizes the log-likelihood function when $\hat\Sigma\not\in S_L$? There is a beautiful geometric answer to this question. This is the unique positive definite matrix $\hat\Sigma_L\in S_L$ such that for all $K\in L$ we have $\langle K,\hat\Sigma_L\rangle=\langle K,\hat\Sigma\rangle$, where the pairing is the trace of the product of matrices. In other words $\hat\Sigma_L-\hat\Sigma\in L^\perp$ \cite[p.~604]{StUh}. For this reason one considers the projection $\pi:\PP(S^2\CC^n)\dashrightarrow\PP(S^2\CC^n/L^\perp)$, where $S^2$ is the second symmetric power, i.e.~$S^2\CC^n$ is identified with the space of symmetric $n\times n$ matrices. The ML-maximization problem then turns out to be related to finding the fiber of the generically finite map $\pi_{|\overline{S_L}}$ over $\pi(\hat\Sigma)$. Indeed, there is a unique positive definite matrix in $S_L$ in the fiber and it is the maximizer of the likelihood function. What about other, possibly complex, points in the fiber? The cardinality of the fiber, for general $\hat \Sigma$, is precisely the ML-degree \cite{amendola2021maximum}. In particular, we see that the ML-degree is the degree of the map $\pi:\PP(\overline{S_L})\dashrightarrow \PP(S^2\CC^n/L^\perp)$. Equivalently, it is the cardinality of the intersection $|(W\cap \PP(\overline{S_L}))\setminus \PP(L^\perp)|$, where $W$ is a general projective subspace of dimension equal to the codimension of $\PP(\overline{S_L})$ that contains $\PP(L^\perp)$. This is always upper bounded by the degree of the model and equality holds when $\PP(L^\perp)\cap \PP(S_L)=\emptyset$. By a theorem of Teissier \cite[II.2.1.3]{T2}, \cite[Lemma 3.1]{jiang2021linear}, this happens when $L$ is general, however often fails for special $L$.

In algebraic statistics both cases are interesting: general and special $L$. Of particular interest are \emph{graphical Gaussian models}, where a graph $G$ with $n$ vertices, indexed by integers $1,\dots,n$, encodes the space $L^G$. Precisely, $L^G$ is the subspace of symmetric matrices that have zeros on off-diagonal entries $(a,b)$ whenever there is no edge between $a$ and $b$ in $G$. Note that we always allow arbitrary diagonal entries, thus one can imagine loops at every vertex of $G$. We emphasise that this is a very different construction of a subspace $L$ of matrices from the one described in the previous section and denoted by $L_G$. Before we had a subspace of diagonal matrices, while now we obtain a space of symmetric matrices that contains all diagonal matrices. We thus have the following interesting invariants:
\begin{itemize}
\item the degree of the model for general $L$ of codimension $a$, which is the degree of the variety $\overline{S_L}$. This degree also equals the ML-degree and is denoted by $\phi(n,a)$,
\item the degree of the model $L^G$,
\item the ML-degree of $L^G$. 
\end{itemize}
The last two of course depend on $G$. It turns out that even for relatively simple $G$, computing those invariants is highly nontrivial.
\begin{conjecture}\cite[7.4]{drton2008lectures} If $G$ is an $n$-cycle then the ML-degree equals:
\[(n-3)\cdot 2^{n-2}+1.\]
\end{conjecture}
The following theorem was proved using the methods described at the end of the article and confirms a conjecture of Sturmfels and Uhler \cite{StUh}.
\begin{theorem}\cite[Theorem 1.4]{dinu2021geometry} If $G$ is an $n$-cycle then the degree of the model equals:
\[\frac{n+2}{4}{\binom{2n}{n}-3\cdot 2^{2n-3}}.\]
\end{theorem}
Recall that the map $\PP(S^2\CC^n)\dashrightarrow \PP(S^2\CC^n)$ that is given by inverting the matrix, is also the gradient of the determinant, i.e.~the map given by all partial derivatives of the determinant.
We have the following commutative diagram \cite[p.6]{dinu2021applications}:
\begin{equation}\label{diagram}
\begin{tikzcd}  
 \PP(S^2\CC^n) \arrow[rr,"\nabla \det",dashed]  & &\PP(S^2\CC^n) \arrow[d,"\pi",dashed]    \\
\PP(L) \arrow[rr,"\nabla(\det_{|\PP(L)})"',dashed] \arrow[hookrightarrow]{u} \arrow[rru,"(\nabla\det)_{|\PP(L)}",dashed] & & \PP(L^*)=\PP(S^2\CC^n/L^\perp)
\end{tikzcd}
\end{equation}
The degree of the model is the last integer in the multidegree of the graph of the diagonal map $(\nabla\det)_{|\PP(L)}$, while the ML-degree is the last integer in the multidegree of graph of the lower map $\nabla(\det_{|\PP(L)})$. 

We would like to point out that there are other interesting questions regarding the geometry of the model, like the generators of the ideal of $S_{L^G}$, however this lies outside the scope of this article. 
\subsection*{Euler characteristic}
The Euler characteristic $\chi(X)$ is one of the most important discrete invariants of a topological space $X$. Assuming we may triangulate the space, we could define it as the number of vertices in the triangulation minus the number of edges plus the number of triangles etc. For complex algebraic constructible sets Euler characteristic is additive:
\begin{equation}\label{Chi+}\chi(X)=\chi(X\setminus Y)+\chi(Y).\end{equation}
One may thus often compute it by breaking an algebraic set into pieces that are easier to understand. Another useful trick is to equip a projective variety $X$ with an action of an algebraic torus $T=(\CC^*)^k$. Then $\chi(X)=\chi(X^T)$, where $X^T$ is the locus of torus fixed points. In particular, if there are finitely many $T$-fixed points, then the Euler characteristic equals the number of those points. For $X=\PP^n$ with the standard multiplicative action of the torus consisting of points with all coordinates nonzero, we get $\chi(\PP^n)=n+1$. 

Consider a projective hypersurface $X=V(f)\subset\PP^n$. We will be mostly interested in the case when $\PP^n$ is a subspace of the space of matrices and $f$ is the restriction of the determinant. The bottom row in Diagram \eqref{diagram} is thus a special case of the map:
\begin{equation*}
\begin{tikzcd}  
 \PP^n \arrow[r,"\nabla f",dashed] &\PP^n .  \\
 \end{tikzcd}
\end{equation*}
Huh \cite[p.912]{Huh1}, following \cite{dimca2003hypersurface}, observed a beautiful relation among the Euler characteristic and the multidegrees $\mu_i$ of the graph of $\nabla f$:
\begin{equation}\label{eq:Euler}
\chi(\PP^n\setminus X)=\sum_i (-1)^i \mu_i.
\end{equation}
A well-known special case is when $X$ is smooth. Then $\nabla f$ is defined everywhere. To compute $\mu_i$ we may use the isomorphism of the domain $\PP^n$ with the graph. We then have to intersect $i$ general hyperplanes, with $n-i$ polynomials, that are general linear combinations of partial derivatives of $f$. These partial derivatives do not have common zeros, thus by a theorem due to Bertini the intersection is smooth, and if $\deg f=d$, we have $\mu_i=(d-1)^{n-i}$ points, which is the expected B{\'e}zout bound. By \eqref{Chi+} and \eqref{eq:Euler}, we obtain:
\[\chi(X)=\chi(\PP^n)-\chi(\PP^n\setminus X)\]
\[=(n+1)-\frac{1-(1-d)^{n+1}}{d}.\]
When $X$ has isolated singularities, then all but one $\mu_i$'s remain the same as in the smooth case. Namely as long as we intersect with at least one hyperplane in the domain, we will not see the singularities and we may still apply Bertini theorem obtaining the B\'ezout bound. However, for $\mu_n$, the derivatives of $f$ intersect in $\mu_n$ many simple points and in the singular points of $X$. Thus, from $(d-1)^n$ we have to subtract the sum of multiplicities of the isolated singularities of $X$. These multiplicities are known as \emph{Milnor numbers} and in our case are nothing else than the dimension, as a vector space, of the algebra of the singular locus. 
We emphasise that Equation \eqref{eq:Euler} works for arbitrary singularities.

This theory may be applied to answer Question 3 from the introduction. Indeed, a matrix with linear entries, may be identified with a space $L$ of matrices. We may compute $\chi(\PP(L)\cap V(\det))$ from the multidegree of the graph of the map $\PP(L)\dashrightarrow\PP(L^*)$ given by the gradient of the determinant composed with the projection from $L^\perp$. 

\begin{example} 
Consider the space $L$ of $4\times 4$ diagonal, traceless matrices, which we obtained from the $4$-cycle. We identify $\PP(L)=\PP^2\subset\PP^3$. On $\PP^3$ the determinant is the product of the four linear forms corresponding to the coordinates. Thus, the hypersurface $X=V(\det)$ is the hyperplane arrangement given by four planes. When restricted to $\PP(L)$ we obtain four lines, which pairwise intersect. The six intersection points correspond to the $6=\binom{4}{2}$, unique up to scalar multiplication, rank two, diagonal, traceless matrices. Denoting a point by $pt$, we have:
\begin{align*}\chi(\PP(L)\setminus X)=\chi(\PP^2)-4\cdot\chi(\PP^1)+6\cdot\chi(pt)\\
=3-4\cdot 2+6=1.
\end{align*}
This, as expected, coincides with the signed sum of multidegrees, computed in Examples \ref{exm:cycle} and \ref{MD_C4}:
\[1-3+3=1.\]
\end{example}
\begin{example}
Let $L$ be the space of $2\times 2$ matrices. We have $\PP(L)=\PP^3$ and the gradient of the determinant is a linear isomorphism. In particular, all $\mu_i$ are equal to one and their signed sum equals $0$. The determinantal hypersurface $X$ is the quadric from Example \ref{ex:quadric}. We obtain:
\begin{align*}
0=\chi(\PP^3\setminus X)=4-\chi(X).
\end{align*}
Thus, $\chi(X)=4$. Indeed, readers with experience in algebraic geometry may recognize that $X$ is isomorphic to $\PP^1\times \PP^1$, via the Segre embedding. 
\end{example}
\subsection*{Quadric hypersurfaces}
Every homogeneous degree two polynomial $f\in \CC[x_1,\dots,x_n]$ is uniquely represented by a symmetric $n\times n$ matrix $M_f$:
\[f(x)=x^t M_f x.\]
The matrix $M_f$ has on the diagonal coefficients of $x_i^2$ in $f$ and on $(i,j)$-th off-diagonal entry, half of the coefficient of $x_ix_j$. Clearly the properties of $f$ and the variety $V(f)$ are related to the properties of $M_f$. We first learned this relation in school, where, for $n=2$ and $f=ax^2+bxy+cy^2$ we compute $\Delta=b^2-4ac$, which is nothing else than $4\det M_f$. We know that $f=0$ has two distinct complex solutions if and only if $\Delta\neq 0$, i.e.~$M_f$ has rank two.

From now on we identify the space of degree two homogeneous polynomials in $n$ variables with $S^2\CC^n$. We will always assume $f\neq 0$. In general, the quadric $V(f)$  is smooth\footnote{Here, we consider $V(f)$ as a scheme, thus if $f$ is not reduced, i.e.~a square of a linear form, we say it is not smooth.} if and only if $\det M_f\neq 0$ and $f$ is a square of a linear form if and only if rank of $M_f$ equals one. Given a smooth quadric $V(f)$ we may consider all hyperplanes $H\subset \PP^{n-1}$ that are tangent to $V(f)$ at some point. Each such $H$ corresponds to a point $P_H\in (\PP^{n-1})^*$ in the dual projective space. It turns out that the locus of all $P_H$ for which $H$ is tangent to $V(f)$, is also a quadratic hypersurface known as the \emph{dual quadric}. The matrix associated to the dual quadric is, up to scaling, $M_f^{-1}$. 

Our aim is to head towards Question 4 from the introduction. First, fix a point $P_1\in \PP^{n-1}$. Which quadrics pass through $P_1$? Note that $f(P_1)=0$ is a \emph{linear} equation in the coefficients of $f$. Thus we obtain a hyperplane $H_{P_1}\subset \PP(S^2\CC^n)$ of polynomials $f$ such that $P_1\in V(f)$. When we have more points we simply have to intersect $H_{P_1}\cap\dots\cap H_{P_k}$ to obtain the locus of quadrics that pass through $P_1,\dots,P_k$. If $P_i$'s are general we see that: 
\begin{itemize}
\item as long as $k<\dim  \PP(S^2\CC^n)=\binom{n+1}{2}-1$ we have infinitely many quadrics,
\item if $k=\binom{n+1}{2}-1$ there is precisely one quadric, up to scaling,
\item  if $k>\binom{n+1}{2}-1$ there are no quadrics
\end{itemize}
that pass through all $P_i$'s. This is slightly less obvious than one may think, as general $P_i$'s do \emph{not} give general hyperplanes $H_{P_i}$. The formal proof could for example rely on the fact that $H_{P}$'s do not have base locus, i.e.~$\bigcap_{P\in\PP^n}H_P=\emptyset$, which is equivalent to the fact that no quadric passes through all points, plus simple linear algebra. Just to point out what may go wrong we present the following example.
\begin{example} We would like to answer the question: how may, up to scaling, homogeneous degree two polynomials $f=ax^2+bxy+cy^2$ 
\begin{itemize}
\item have a double root and 
\item vanish at a given point $P\in \PP^1$? 
\end{itemize}
We consider the $\PP^2$ of all polynomials. In this $\PP^2$ the locus of polynomials that have a double root is the quadratic hypersurface $\Delta=0$. We intersect this hypersurface with the line $H_P$. We expect two solutions. This is clearly wrong! If we fix a root $P$ and require that a degree two polynomial has a double root, then $P$ must be the double root.
Every child knows that there is only one, up to scaling, such a polynomial. The geometry here is that the line $H_P$ will be always tangent to $V(\Delta)$ and thus it will intersect it in a single point. 
\end{example}
How about the locus of quadrics that are tangent to a given hyperplane $H\subset\PP^{n-1}$? We already know that a smooth quadric represented by $M_f$ is tangent to $H$ if and only if the dual quadric represented by $M_f^{-1}$ passes through the point $P_H\in (\PP^{n-1})^*$. As we are working up to scaling we may replace $M_f^{-1}$ by the adjugate matrix $M_f^{adj}$. The entries of $M_f^{adj}$ are degree $(n-1)$ polynomials in the original coordinates, i.e.~the entries of $M_f$. The condition $(P_H)^t M_f^{adj} (P_H)=0$ is thus a degree $(n-1)$ polynomial defining a hypersurface $T_H\subset \PP(S^2\CC^n)$. Very explicitly this hypersurface is a linear combination of $(n-1)\times (n-1)$ minors of $M_f$. It thus makes sense to ask: how many quadrics are tangent to a given hyperplane and pass through $\binom{n+1}{2}-2$ fixed general points?
Geometrically we intersect $T_H$ with hyperplanes and get $\deg T_H=n-1$ many points\footnote{For experts: as our choices are not entirely generic, we should refer here to Kleiman's transversality theorem, assuring that we indeed obtain $(n-1)$ many points corresponding to smooth quadrics.}.

We could be tempted to continue this game. If we consider $a$ general points and $\binom{n+1}{2}-a-1$ general hyperplanes, how many quadrics pass through the given points and are tangent to the given hyerplanes? Geometrically we intersect $a$ hyperplanes and $\binom{n+1}{2}-a-1$ hypersurfaces of degree $(n-1)$. Thus we expect that the answer is equal to the B{\'e}zout bound: $(n-1)^{\binom{n+1}{2}-a-1}$. This, however, is wrong in general!

The easiest way to see this is to consider $n=3$ and five general lines in $\PP^2$. If we ask for quadrics that are tangent to all five, we may equivalently ask for dual quadrics passing through five points. But here we already know the answer is one, not $32$. What goes wrong? We recall that our hypersurfaces $T_H$ are defined by linear combinations of minors of size $(n-1)$. In particular, each $T_H$ contains all matrices of rank at most $(n-2)$. This is a very different situation than for the hyperplanes $H_P$, which had no base locus. The codimension of the locus of symmetric $n\times n$ matrices of rank $(n-2)$ equals three. Thus, if we take $b\geq 3$ general hyperplanes and intersect the corresponding hypersurfaces $T_H$ we will get:
\begin{itemize}
\item a big codimension three component of matrices of rank at most $(n-2)$ and
\item a small codimension $b$ component.
\end{itemize}
The meaningful geometric counting problem is to ask for \emph{smooth} quadrics that pass through the given points and are tangent to the given hyperplanes. Thus, we would like to intersect the small codimension $b$ component with hyperplanes $H_{P_i}$ and count the number of points. However, the contribution of the big component makes such computations quite hard.

Let us state the enumerative problem in a way more suitable for this article. By fixing points $P_1,\dots, P_a$, we fix a linear space $L\subset \PP(S^2\CC^n)$ of quadrics that pass through those points. We may now consider the rational map:
\[\PP(S^2\CC^n)\dashrightarrow\PP(S^2\CC^n)^*\]
that is the gradient of the determinant, or equivalently, taking adjugate or inverse of a matrix\footnote{In the projective setting all these maps are the same on the Zariski open set of invertible matrices. However, taking adjugate matrix is well-defined also for matrices of rank $n-1$.}. We restrict the map to $\PP(L)$, obtaining the diagonal map in Diagram \eqref{diagram}. Each tangency to a hyperplane condition is in fact intersection with a hyperplane in the dual space $\PP(S^2\CC^n)^*$. We thus see that the number of quadrics that pass through the given points and are tangent to the correct number of general hyperplanes is in fact one entry in the multidegree of the graph of $(\nabla\det)_{|\PP(L)}$. Note that here the problem of base locus and low rank matrices disappears. Indeed in the product of projective spaces both: hyperplanes in the domain and hyperplanes in the codomain do not have the base locus. Also as the graph is by definition the closure of the locus corresponding to full rank matrices, we do not have the additional large component. Note that for general points $P_1,\dots,P_a$ we get:

\noindent\emph{the number of quadrics that pass through all $P_i$'s and are tangent to general $\binom{n+1}{2}-a-1$ hyperplanes equals precisely the degree $\phi(n,a)$ of the general linear concentration model.}

In many cases in enumerative problems the way to obtain the correct answer in a nice way is to change the ambient space where intersection is performed. This usually requires very good, clever ideas. Passing to the graph of the map, instead of $\PP^n$ or $\PP(L)$ seems to be one. However, it is only the first step, as such graph is in general not smooth.
Fortunately, in the cases interesting to us, great mathematicians before us already had the right ideas. We may present these in the next subsection.
\subsection*{Complete varieties}
In a quite non-standard way, we do \emph{not} refer to complete varieties, as a synonym for proper. The complete varieties as described below are very special projective varieties, that may be regarded as particularly nice compactifications of the locus of nondegenerate matrices.  

Let $W$ be the space of diagonal or symmetric or general $n\times n$ matrices. Let $W^\circ$ be the subset of full rank matrices. We could see $\PP(W^\circ)\subset\PP(W)$ as a natural compactification. However, as we have seen in the previous section the low rank matrices in $\PP(W)$ often turn out to be problematic. There exists a well-known procedure in algebraic geometry of replacing a small set by a divisor, i.e.~a codimension one set, known as the \emph{blow-up}. Say we have a subset $S\subset\PP^n$, where $S=V(g_0,\dots,g_s)$ and all $g_i$'s are polynomials of fixed degree $d$. For our purposes, we define the blow-up of $\PP^n$ at $S$ as the graph of the rational map given by $g_i$'s:
\[\PP^n\dashrightarrow \PP^s.\]
This procedure is particularly nice when $S$ is smooth. In our case $\PP^n=\PP(W)$ and $S$ is the set of rank at most $n-2$ matrices. Thus, the graph we considered, is the blow-up construction, where $g_i$'s are the minors. However, the set $S$ is not smooth, it is singular along matrices of rank at most $n-3$, which further is a set singular along matrices of rank at most $n-4$ etc.~until rank one matrices, which is a smooth locus in $\PP(W)$. 

The idea is to first blow-up rank one matrices, then rank two matrices, etc.~until finally we blow-up rank $n-2$ matrices. This may be realized at once as follows. Consider the rational map:
\[\psi:\PP(W)\dashrightarrow \PP(W_2)\times\dots\times \PP(W_{n-1})\]
where the map $\PP(W)\dashrightarrow \PP(W_i)$ is given by $i\times i$ minors. The dimension of $\PP(W_i)$ depends on the case we are in, e.g.~for diagonal matrices $\dim W_i=\binom{n}{i}$, while for general matrices $\dim W_i=\binom{n}{i}^2$, however otherwise the construction remains the same. The (closure of the) graph of $\psi$ is known as:
\begin{itemize}
\item permutohedral variety in the diagonal case,
\item variety of complete quadrics in the symmetric case,
\item variety of complete collineations in the case of general square matrices\footnote{One may also define this variety for rectangular matrices, however here we do not pursue this direction.}.
\end{itemize} 
It is equal to the iterative blow-up of rank $i$ matrices for $i=1,\dots,n-2$ anticipated above. This variety is smooth. By projecting to $\PP(W)\times\PP(W_{n-1})$ it becomes a resolution of singularities of the graph of the map inverting the matrix. In this article, to emphasise the common features of the above varieties, we refer to all of them as \emph{complete varieties}.

As an example let us show how one may think about points of the variety $X$ of complete quadrics. First, there is an open locus, corresponding to quadrics of rank at least $(n-1)$, which is isomorphic to an open subset of $\PP(W)$, via the projection $X\rightarrow \PP(W)$. However, there are many points $x\in X$ mapping to a given quadric $Q\in \PP(W)$ if the rank of $Q$ is at most $n-2$. Fix such a quadric $Q\in \PP(W)$. Let $V\subset \CC^n$ be the image of the symmetric matrix $A_Q$. It turns out that the fiber over $Q$ is the variety of complete quadrics over the vector space $\CC^n/V$. We see that the points of the variety of complete quadrics may be identified with the data consisting of:
\begin{itemize}
\item  flags $V_1\subset\dots\subset V_k=\CC^n$ and
\item full rank quadrics on each $V_{i+1}/V_i$, i.e.~elements of $\PP(S^2 (V_{i+1}/V_i)^\circ)$.
\end{itemize}
A particularly interesting case is when the flag is full, as then the data of quadrics is trivial --- there is only one up to scaling nonzero, degree two homogeneous polynomial in one variable. Thus we see that the variety of complete quadrics naturally contains the full flag variety $F$. For experts: the inclusion $F\rightarrow X$ induces the inclusion of Picard groups $\Pic(X)\rightarrow\Pic(F)$. The right hand side is understood in terms of the root system of type $A$. This allows us to have a good understanding of $\Pic(X)$. 


\section{Intersection theory}

Classical algebraic intersection theory associates to an algebraic variety $X$, the graded Chow ring $\CH(X)$. We have already seen this on the example of $X=\PP^n$. In general, the degree $k$ part of $\CH(X)$ is spanned by classes of codimension $k$ subvarieties, modulo rational equivalence. For example for $k=1$ we obtain the divisor class group, which for smooth $X$ coincides with the Picard group $\Pic(X)$. The multiplication in $\CH(X)$ corresponds to intersection. Precisely, if $Y,Z\subset X$ intersect transversally, then $[Y]\cdot[Z]=[Y\cap Z]$. For all the varieties that we consider in this article the Chow ring is isomorphic to the cohomology ring. We refer the reader interested in intersection theory to the book \cite{eisenbud20163264}.

\subsection*{Volume polynomials}
Given a $k$-dimensional subvariety $Y\subset X$ we may define a polynomial function on the Picard group $\Pic(X)$:
\[\Pic(X)\ni D\mapsto \deg ([Y][D]^k)\in \ZZ,\]
where $\deg$ is the degree function that to a zero-dimensional (class of a) scheme associates its degree. If we prefer to work in coordinates and obtain a polynomial we fix divisors $D_1,\dots,D_n\in\Pic(X)$. We define the \emph{volume polynomial} of $Y$ on $\QQ^n$ by:
\[(t_1,\dots,t_n)\mapsto \deg\left((\sum t_i[D_i])^k[Y]\right).\]

A divisor $D$ is called nef (numerically effective) if for every curve $C\subset X$ the number $\deg([C]\cdot[D])$ is nonnegative. A particularly nice case is that of a volume polynomial when $D_1,\dots,D_n$ are nef. Then the volume polynomial is \emph{Lorentzian} \cite[Theorem 4.6]{branden2020lorentzian}, exhibiting a lot of nice properties.

In case of complete varieties there is a distinguished set of nef divisors. By definition, the complete variety is a subvariety of the product $\PP(W)\times\PP(W_2)\times\dots\times\PP(W_{n-1})$, the hyperplane $H_i\subset \PP(W_i)$ (times the remaining $\PP(W_j)$'s) gives a nef divisor $L_i$ on $X$. 

For example, when $X$ is the variety of complete quadrics, the $L_i$'s are the extremal rays generating the nef cone --- all nef divisors are nonnegative linear combinations of $L_i$'s. Via the containment of $\Pic(X)$ in the root system of type $A$ the $L_i$'s correspond to (twice) the fundamental roots. 
For the permutohedral variety the $L_i$'s do not generate the nef cone, however they generate the $S_n$ (i.e.~permutation) invariant part. 

\underline{To sum up}: consider an $a$-dimensional linear space $L$ of $n\times n$ matrices. Assume $L$ contains an invertible matrix. The closure of invertible matrices in $L$ in a complete variety $X$ (which depends on the type of matrices $L$ consists of) gives a subvariety $Y_L$. We obtain \emph{the Lorentzian polynomial} given by $\deg(\sum t_iL_i)^{a-1}[Y_L]$. In \cite{conner2021characteristic} the coefficients of this polynomial were introduced and called the \emph{characteristic numbers}.

By setting $t_2=\dots=t_{n-1}=0$ we recover the \emph{chromatic polynomial} of a tensor, defined also in  \cite{conner2021characteristic}. Its coefficients (up to binomial factors) are precisely the multidegrees of the graph of the map inverting matrices\footnote{Equivalently gradient of the determinant} from $L$. 

In a very related construction, replacing the restriction to $L$ of the gradient of the determinant, by the gradient of the restriction to $L$, we obtain the \emph{relative chromatic polynomial}. These two coincide e.g.~when $L$ consists of diagonal matrices, or when $L$ is general.  These two invariants seem very important and, in particular, provide answers to all the questions from the introduction.

\subsection*{Methods in a nutshell}
Here, we would like to briefly mention a few methods that rely on complete varieties and allow to compute the intersection numbers. For more details we refer to \cite{manivel2020complete, michalek2021maximum, dinu2021applications}. 
For simplicity, consider the variety $X$ of compete quadrics. Say we want to compute $[L_1]^a[L_{n-1}]^{\binom{n+1}{2}-1-a}$. The first trick is due to Schubert. One notices that on $X$ we have another set of divisors. Indeed, let $S_i$ be the exceptional divisor coming from the blow-up of rank $i$ matrices. It turns out that $S_i$ correspond to (twice) the simple positive roots. We thus have relations in $\Pic(X)$ that allow us to translate from $L_i$'s to $S_i$'s. We obtain:
\[[L_1]=\frac{1}{n}\sum_{j=1}^{n-1}(n-j)[S_{j}].\]
This reduces the problem to computing all $[S_j][L_1]^{a-1}[L_{n-1}]^{\binom{n+1}{2}-1-a}$. At first this looks like a more complicated problem. The trick is to realize that one can replace $[S_j]$ by a simpler variety. First, we have a rational map $f_j:S_j\dashrightarrow G(j,n)$, from $S_j$ to the Grassmannian parameterizing $j$-dimensional subspaces of an $n$-dimensional vector space.
To a point $s\in S_j$ mapping to a rank $j$ symmetric matrix $M_s$ the map $f_j$ associates the image of $M_s$. Fixing a point $V\in G(j,n)$ we see that the fiber of $f_j$ is birational with: $\PP(S^2V)$ for the choice of $M_s$ times $\PP(S^2(\CC^n/V)^*)$, as after fixing $M_s$ we said that the fiber of $X$ over $M_s$ is isomorphic to the variety of complete quadrics on $\CC^n/V$. Hence, $S_j$ is birational to the product bundle $\PP(S^2\mathcal{U})\times \PP(S^2\mathcal{Q}^*)$ over $G(j,n)$, where $\mathcal{U}$ and $\mathcal{Q}$ are respectively the universal and quotient bundle. This birational morphism is good enough to switch from intersecting with $[S_j]$ on $X$ to intersecting on the product bundle. The classes $[L_1]$ and $[L_{n-1}]$ also translate nicely and may be expressed in terms of characteristic classes of  $S^2\mathcal{U}$ and $S^2\mathcal{Q}$. These have been investigated in detail, cf.~\cite{pragacz1988enumerative} and references therein. 
\section{Future}
\begin{enumerate}
\item When $L$ is general, in the diagonal case, the associated matroid is uniform and characteristic numbers are simply binomial coefficients $\binom{n}{a}$. For fixed $a$, this is a polynomial in $n$ with very nice reciprocity properties. In the symmetric case we obtain $\phi(n,a)$, which is also a polynomial in $n$. However, here the reciprocity results are much less understood, with first, very recent results presented in \cite{Galazka}. For positive $n$ and $a$ the function $\phi(n,a)$ counts the number of quadrics satisfying given conditions. Is there a nice interpretation for fixed positive $a$ and negative $n$?
\item As to  a $1$-generic tensor we associate a Lorentzian polynomial, we obtain a finite subdivision of the space of $1$-generic tensors. What is this subdivision? What are closures of different loci? This study was initiated in \cite{Hanieh}. One may also ask for other algebraic invariants of $L^{-1}$. 
\item One can check that $\deg L^{-1}$ is semicontinuous. In particular, we may define it for any variety in a tensor space, like tensors of bounded border rank. It seems a highly nontrivial task to decide if the invariants coming from Lorentzian polynomials may distinguish cactus border rank from border rank.
\item In the nondiagonal cases, e.g.~when $X$ is the variety of complete quadrics, the Chow ring $\CH(X)$ does not have to be generated in degree one. Thus one obtains more invariants by intersecting the class of $Y_L$ with other (higher degree) classes in the Chow ring. These invariants, as far as we know, have not been studied. 
\item One may enhance the intersection theory using group actions. In the diagonal case the group is the torus, which is heavily exploited in many different settings \cite{berget2021tautological}. In particular, the numbers we obtain are in fact formal combinations of dimensions of representations of the torus. This is true even for the easiest, general case, when the answer is $\binom{n}{a}$. 
A more detailed study of the associated representations is still very much work in progress. The next step will be to go to the case of other complete varieties and study associated invariants.
\item In combinatorics, there is a beautiful theory of duality, which in our setting corresponds to replacing $L$ by $L^\perp$. For tensors this is essentially the castling transform.
At this point, it is not obvious how the invariants behave under this duality for tensors. Also in combinatorics deletion and contraction are central in the study of chromatic polynomials. In principle, one may define contraction with respect to a rank one matrix and deletion with respect to a rank one matrix in the dual space. This would satisfy the well-known fact that the dual of deletion is contraction of the dual. However, relations to Lorentzian polynomials at this point are not clear.  
\item What is the correct way to generalise the stellahedron from \cite{eur2022stellahedral} to nondiagonal cases?
\item The \emph{Bodensee program} \cite{dinu2021applications} studies how discrete invariants, like $\phi(n,a)$, change, as the variety, here $n$, changes. We already mentioned the polynomiality result of $\phi(n,a)$, however the coefficients of this polynomial remain quite mysterious. Their absolute values seem to be log-concave. Is there a hidden cohomology theory, possibly on an infinite dimensional variety, explaining this phenomenon?
\item The cases described in this article correspond to \emph{representable} matroids. However, the cohomology ring of the permutohedral variety or stellahedral variety allow also study of general matroids. Is there a way to upgrade this theory, say replacing permutohedral variety by the variety of complete quadrics, to get a new theory of anabelian matroids?
\item A finite $\CC$-algebra $A$ naturally defines a bilinear map $A\times A\rightarrow A$, given by multiplication. In particular, the algebra $A$ defines a tensor. Hence, the constructions from this article associate a Lorentzian polynomial to such an algebra. What kind of properties of $A$ are encoded by this polynomial? What is the induced Lorentzian subdivision of the Hilbert scheme? 
\end{enumerate}


\bibliography{bibML}
\bibliographystyle{plain}

\end{document}